\documentclass{amsart}[12pt]
\usepackage{amssymb,amsfonts,amsmath,latexsym,mathtools,enumerate}
\usepackage{xcolor}
\usepackage{url}
\newtheorem{theorem}{Theorem}
\newtheorem{lemma}[theorem]{Lemma}

\newtheorem{proposition}[theorem]{Proposition}
\newtheorem{remark}[theorem]{Remark}
\newtheorem{corollary}[theorem]{Corollary}

\newcommand\N{{\mathbb N}}
\newcommand\R{{\mathbb R}}
\newcommand\C{{\mathbb C}}

\renewcommand\S{{\mathbb S}}

\newcommand\balpha{\boldsymbol\alpha}
\newcommand\bbeta{\boldsymbol\beta}

\newcommand\ind{{\mathbf 1}}
\renewcommand\k{{\mathbf k}}

\newcommand\K{{\mathbf K}}
\newcommand\E{{\mathbb E}}

\newcommand\eps{\varepsilon}
\newcommand{\Cov}{\mbox{\rm Cov}}
\newcommand{\Var}{\mbox{\rm Var}}
\newcommand{\vol}{\mbox{\rm vol}}
\newcommand{\sinc}{\mbox{\rm sinc}}
\newcommand{\diag}{\mbox{\rm diag}}

\renewcommand{\>}{\right\rangle}

\begin{document}

\title{On $3$-dimensional Berry's model}
\date{\today}

\author{F. Dalmao}
\address{
DMEL, 
Universidad de la Rep\'{u}blica, Salto, Uruguay.}
\email{fdalmao@unorte.edu.uy.}

\author{A. Estrade}
\address{Laboratoire MAP5, Universit\'e Paris Descartes, Paris, France.} 
\email{anne.estrade@parisdescartes.fr.}

\author{J. R. Le\'{o}n}
\address{IMERL, Universidad de la Rep\'{u}blica, Montevideo, Uruguay and Escuela de Matem\'{a}tica. Facultad de Ciencias. 
Universidad Central de Venezuela, Caracas, Venezuela.} 
\email{rlramos@fing.edu.uy }
\maketitle

\begin{abstract}
This work aims to study the dislocation or nodal lines of 3D Berry's random wave model. 
Their expected length is computed both in the isotropic and anisotropic cases, being them compared.  
Afterwards, in the isotropic case the asymptotic variance and distribution of the length are obtained as the domain grows to the whole space. Under some integrability condition on the covariance function, a central limit theorem is established. 
The study includes the Berry's monochromatic random waves, the Bargmann-Fock model and the Black-Body radiation as well as a power law model that exhibits an unusual asymptotic behaviour and yields a non-central limit theorem.
\end{abstract}

\noindent 
{\bf AMS classification:} 60G60; 60G15, 60F05, 60D05.\\
{\bf Keywords:} Random waves, nodal statistics, central limit theorem, non-central limit theorem.


\section{Introduction}
In the last few years nodal or dislocation (i.e: zero) sets of several classes of random waves have 
received a lot of attention from Number Theory, Topological Analysis, Differential Geometry, Probability Theory, etc. 
While studying the random billiards, Berry \cite{Be02} argued that in the microscopic scale several models 
as arithmetic random waves on the torus or spherical harmonics, 
although they verify some boundary conditions, converge towards an universal Gaussian model, 
which is called Berry's random waves model.  
Canzani and Hanin \cite{CH16} studied the universality phenomenon in general Riemannian manifolds. 
The reader can find results on arithmetic random waves defined on the flat torus \cite{C17,DNPR17}
and on random spherical harmonics in \cite{CM19,fcmt,mr} and references therein, see also \cite{r18} for a survey on both subjects.
The nodal sets of Berry's planar random waves, {\em i.e.} the random eigenfunctions of the 2D Euclidean Laplacian operator, have been studied in \cite{NPR17} where Central Limit Theorems are obtained for the nodal length in the real case and for the number of phase singularities in the complex case. 
Whereas all the previous references are concerned with 2-dimensional isotropic random fields, one can also find studies in more general frameworks. In \cite{EF18} anisotropic random waves are considered in any dimension. In \cite{KV18,Mu17}, similar central limit results are obtained for any Minkowski functional of excursion sets in the general framework of stationary Gaussian fields whose covariance function is fast enough decreasing at infinity.

Our motivation mostly comes from Berry and Dennis seminal paper \cite{BD00} and from Dennis \cite{De07} 
where the authors show how the expectation and the second moment of certain functionals of the nodal sets can be computed. 
The main tools are the different forms of the Kac-Rice formulas 
(see \cite{aw} and the references therein).  
Moreover, in \cite{De07} (where a more formal approach from the mathematical point of view is presented) 
a variety of problems which are in close relation with the computation of the measure of the zero set of random waves  are exhibited. Also, the two points correlation is introduced defining it as a second order Rice's function.
\\

In the present paper we study complex-valued 3-dimensional Berry's random waves models with a focus on the length of the dislocation or nodal lines. We obtain the expected length in a very general framework which includes anisotropy. 
In order to study the asymptotic variance and the limit distribution, we let the domain increase to the whole space. It can be shown that this is equivalent to consider a fixed domain and taking the high energy limit, see \cite{CH16,NPR17}. 
We establish the order of magnitude of the limit variance and the asymptotic normality in a more restrictive framework 
including Berry's monochromatic random waves, Black-Body radiation and Bargmann-Fock waves. 
We also include a power law model which has an asymptotic variance of different order and 
which presents a non-Gaussian limit distribution yielding a non-central limit theorem. 
It is worth to say that in the monochromatic case we find that the behaviour of the asymptotic nodal length variance 
differs from the two dimensional case. Indeed, in that case, the asymptotic variance scales proportionally to the volume, 
without any logarithmic term as in 2D. Nevertheless, the second chaos vanishes as in 2D.
\\

The paper is organized as follows. Section \ref{sec:model} presents the model as well as our main results, namely Theorem \ref{TheoLength} and Theorem \ref{TheoVar}. Section \ref{sec:exp} is devoted to the study of the first moment of the dislocation length; in particular, it contains the proof of Theorem \ref{TheoLength}. A special section (Section \ref{sec:H}) is devoted to the It\^o-Wiener's chaotic decomposition given by Hermite expansion. In Section \ref{sec:clt}, the asymptotic results and proofs are gathered. First, we prove Theorem \ref{TheoVar}. Secondly, we exhibit examples of three dimensional isotropic random waves models with 
various asymptotic behaviours. The two specific cases of Berry's monochromatic random waves and power law model are studied in Sections \ref{sec:Berry} and \ref{sec:power} respectively.

\section{Model and main results} \label{sec:model}
Consider a $3$-dimensional Berry's random waves model $\psi:\R^3\to\C$ given by 
\begin{equation}\label{e:psi}
\psi(x)=\int_{\R^3}\exp(i \langle \mathbf{k},x\rangle)\,
\frac{dW_\Pi(\mathbf{k})}{|\mathbf{k}|},\quad x\in \R^3,
\end{equation}
where $\langle\cdot,\cdot\rangle$ and $|\cdot|$ stand for the usual inner product 
and $2$-norm in $\R^3$ respectively. 
Besides, $W_\Pi$ is a (complete) complex-valued Gaussian random measure 
on $\R^3$ with (real) control measure $\Pi$, {\em i.e.} $\Pi$ is a positive non-atomic measure on $\R^3$ satisfying 
\begin{equation}\label{eq.control}
\E\left(\int_{A}\frac{dW_\Pi(\k)}{|\k|}\,\overline{\int_{B}\frac{dW_\Pi(\k)}{|\k|}}\right)
=2\int_{A\cap B}\frac{\Pi(d\k)}{|\k|^2},
\end{equation}
for any Borel sets $A,B$ in $\R^3$. 
We further assume that  
\begin{equation*}
\int_{\R^3}\frac{\Pi(d\k)}{|\k|^2}=1,
\end{equation*}
that $\Pi(\R^3)<\infty$ 
and that 
$\Pi(-A)=\Pi(A)$ for any Borel set $A\subset \R^3$.

Actually, if $W_\Pi=W_\Pi^1+iW_\Pi^2$ with real independent $W_\Pi^j$ ($j=1,2$) then 
\eqref{eq.control} holds for $W_\Pi^j$ ($j=1,2$) without $2$ in the right-hand side factor. \\

As a consequence, the random field $\psi$ is Gaussian, stationary, centered but not necessarily isotropic. 
We denote by $\xi$ and $\eta$ the real and imaginary parts of $\psi$, that is $\psi=\xi+i\eta$. 
The random fields $\xi$ and $\eta$ are independent and identically distributed 
with common covariance function prescribed by
\begin{eqnarray}\label{e:r}
r(x)&:=&\E(\xi(0)\xi(x))=\E(\eta(0)\eta(x)),\quad x\in \R^3 \notag\\
&=&\int_{\R^3}\exp(i\langle \mathbf{k},x\rangle)\,\frac{\Pi(d\k)}{|\k|^2}.
\end{eqnarray}
Note that the normalization $\int_{\R^3}\frac{\Pi(d\k)}{|\k|^2}=1$ yields $r(0)=1$ and that
$$\E(\psi(0)\overline{\psi}(x))=2r(x).$$ 
Furthermore, the condition $\Pi(\R^3)<\infty$ implies that $\psi,\xi,\eta$ are almost surely $C^2$.

\medskip

Using the vocabulary introduced in \cite{EF18}, $\xi$ and $\eta$ are 
random waves whose associated random wavevector admits $\frac{\Pi(d\k)}{|\k|^2}$ as distribution. 
In what follows, we will call $\Pi$ the \textit{power spectrum} although this word is usually reserved 
to the isotropic framework. 
Indeed, the random wave $\psi$ can be isotropic or not 
according to the fact that the covariance function $r(x)$ only depends on $|x|$ or not, 
which only depends on the choice of $\Pi$.

\medskip

Let us look at the model in the isotropic case. 
We write $\k=\rho\mathbf{u}$ with $\rho>0$ and 
$\mathbf{u}\in \S^2$, being $\S^2$ the unitary sphere in $\R^3$. 
We consider the case where the image of measure $\frac{\Pi(d\k)}{|\k|^2}$ through the change of variables $\k\mapsto (\rho,u)\in \R^+\times \S^2$ writes out 
\begin{equation}\label{eq:Pirad}
\Pi^{rad}(d\rho)\otimes d\sigma(\mathbf{u}),
\end{equation}
for some measure $\Pi^{rad}$ defined on $\R^+$ 
and where $d\sigma$ stands for the surface measure on $\S^2$. 
The normalization on the power spectrum imposes that $\Pi^{rad}(\R^+)=\frac1{4\pi}$.  
The covariance function is then given by
\begin{eqnarray}\label{e:riso}
r(x)&=&\int_{\R^+}\,\left( \int_{\S^2}\exp(i\rho |x|\langle 
\mathbf{u},\mathbf{e}\rangle)\,d\sigma(\mathbf{u}) 
\right)\,\Pi^{rad}(d\rho)\notag\\
&=&4\pi\,\int_{\R^+}\frac{\sin(\rho|x|)}{\rho|x|}\,\Pi^{rad}(d\rho),
\end{eqnarray}
being $\mathbf{e}$ a fixed point in $\S^2$. 
 
In view of \eqref{e:riso}, we recognize the covariance function involved in Berry and Dennis model \cite{BD00}. 
In particular, if $\Pi^{rad}$ is a Dirac measure at some point $\kappa$ in $(0,\infty)$, we recover Berry's monochromatic random wave. Note that our normalization on $\Pi^{rad}$ differs from (3.11) in \cite{BD00}. 
More examples are gathered in Section \ref{sec:clt}.

\medskip


Let us turn to the main purpose of the paper: the study of the length of the dislocation lines $\{x\in \R^3\,:\,|\psi(x)|=0\},$ 
which have Hausdorff dimension one. For any bounded domain $Q$ in $\R^3$, we introduce
$$
{\mathcal Z}(Q)=\{x\in Q\,:\,|\psi(x)|=0\},\quad
\ell ({\mathcal Z}(Q))={\rm length}({\mathcal Z}(Q)).
$$

We now present our main results. 
The first one deals with the expectation of $\ell({\mathcal Z}(Q))$ for arbitrary fixed $Q$.
\begin{theorem} \label{TheoLength}
Let $\psi$ be defined as in \eqref{e:psi} and assume that $\psi'(0)$ is non degenerated. 
Let $\lambda_i$, $i=1,2,3$ be the eigenvalues of the covariance matrix $-r''(0)$ and $D=\diag(\lambda_1,\lambda_2,\lambda_3)$.
Hence, 
\[
\E(\ell ({\mathcal Z}(Q)))=\frac{\sqrt{\lambda_1\lambda_2\lambda_3}}{2\pi}
\E|D^{-\frac12}(N\wedge N')|\,\vol(Q),
\]
being $(N,N')$ a standard normal random vector in $\R^6$ and $\wedge$ the usual cross product of vectors in $\R^3$.
\end{theorem}

Next, we specialize this result to the isotropic case and compare it with the \emph{almost} isotropic case.  
\begin{corollary}\label{cor:length}
In the same conditions as above, 
\begin{enumerate}
\item [(i)] if $\lambda_i=\lambda=-r''_{11}(0)$, $i=1,2,3$, we have 
\[
\E(\ell ({\mathcal Z}(Q)))=\frac{\lambda}{\pi}\,\vol(Q)~;
\]
\item [(ii)] for $\lambda>0$ fixed, as $\max_i|\lambda_i-\lambda|\to0$, 
we have the following expansion
\begin{multline*}
\phantom{space}\E(\ell({\mathcal Z}(Q)))=\frac{\lambda}{\pi}\,\vol(Q)\left(
1+(-1+\frac{2}{3}\sqrt{\lambda})\sum^3_{i=1}(\lambda_i-\lambda)\right)+O(\max_i|\lambda_i-\lambda|^2).
\end{multline*}
\end{enumerate} 
\end{corollary}

The proofs of Theorem \ref{TheoLength} and Corollary \ref{cor:length} are postponed to Section \ref{sec:exp}. 
The first item in the corollary is coherent with (3.14) in \cite{BD00} taking into account that $\lambda=k_3/3$ in Berry and Dennis notation. Besides, we have
\[
\E(\ell({\mathcal Z}(Q)))=\frac{\E|\xi'(0)\wedge\eta'(0)|}{2\pi}\,\vol(Q),
\] 
where 
$\xi'(0)\wedge\eta'(0)$ is the so-called {\it vorticity}, see  (2.2) in \cite{BD00}.
~\\

We now move to the asymptotic behaviour of the variance and distribution of $\ell({\mathcal Z}(Q))$ as $Q$ 
grows up to $\R^3$. We restrict ourselves to the isotropic case, assuming moreover that the radial component $\Pi^{rad}$ in \eqref{eq:Pirad} admits a density and that the covariance function is square integrable.


Let
\begin{equation} \label{defR}
R(x)=\max\left\{|r(x)|,|r'_i(x)|,|r''_{ij}(x)|\,:\, 1\le i,j\le 3\right\},~x\in \R^3. 
\end{equation}

\begin{theorem} \label{TheoVar}
Let $\psi$ be an isotropic Berry's random wave defined as in \eqref{e:psi} and \eqref{eq:Pirad} such that $\Pi^{rad}$ admits a density with respect to Lebesgue measure. 
Assume that $R(x)\mathop{\to}0$ whenever $|x|\to\infty$ and that $R\in L^2(\R^3)$. Finally, let $Q_n=[-n,n]^3$. 
Hence, 
\begin{enumerate}
\item[(i)] there exists $0<V<\infty$ such that
\begin{equation*}
 \lim_{n\to \infty}\frac{\Var(\ell(\mathcal{Z}(Q_n)))}{\vol(Q_n)}=V;
\end{equation*}
\item[(ii)] as $n\to\infty$, the distribution of  
\[ \frac{\ell({\mathcal{Z}(Q_n)})-\E(\ell({\mathcal{Z}(Q_n)}))}{\vol(Q_n)^{1/2}} \]
converges towards the centered normal distribution with variance $V$.
\end{enumerate}
\end{theorem}
~\\
The proof of Theorem \ref{TheoVar} can be found in Section \ref{sec:L2case}. \\
Performing the isotropic space scaling $x\mapsto \kappa x$ in $\R^3$ for some $\kappa>0$ yields the next remark. 
\begin{remark} \label{rem:highNRJ}
If $\psi$ is an in Theorem \ref{TheoVar} and if $\psi_\kappa$ is defined as $\psi_\kappa=\psi(\kappa\,\cdot)$ then, the distribution of 
\[\frac{\mbox{length}(\psi_\kappa^{-1}(0)\cap [-1,1]^3)-\E(\mbox{length}(\psi_\kappa^{-1}(0)\cap [-1,1]^3))}{\kappa^{1/2}} \]
converges as $\kappa$ tends to $+\infty$ towards a centered normal distribution with some variance $V \in (0,\infty)$.
\end{remark}
One can see this asymptotics either as an {\em infill statistics} statement since the performed scaling is nothing but a zooming (see \cite{CH16}), or as a {\em high energy} statement (see \cite{NPR17}) since the second spectral moment $\lambda_\kappa$ of $\psi_\kappa$ is such that $\lambda_\kappa=\kappa^2\lambda$ and hence tends to $+\infty$.

\medskip

\section{Expected nodal length} \label{sec:exp}
In this section we compute the mean length of the dislocation lines and prove Theorem \ref{TheoLength} 
and Corollary \ref{cor:length}. 

\medskip

We need some further notations.  For any $x\in \R^3$, let $Z(x)=(\xi'(x),\eta'(x))$ where
$Z(x)$ is sometimes considered as a vector in $\R^6$ and sometimes as a $2\times 3$ matrix. We also denote
\[
\det{}^{\bot}Z(x)=\det{}^{\bot}\left(\begin{array}{ccc}\xi'_1(x)&\xi'_2(x)&\xi'_3(x)\\ \eta'_1(x)&\eta'_2(x)&\eta'_3(x)\end{array}\right),
\]
where for any real matrix $M$ $\det{}^{\bot}M$ stands for $\det(MM^{\top})$.
Routine computation shows that
\begin{equation}\label{e:cross}
\det{}^{\bot}Z(x)=|\xi'(x) \wedge \eta'(x)|^2.
\end{equation}
This equality is a particular case of the well known Binet-Cauchy formula. 
\\~\\
The expectation of $\ell({\mathcal Z}(Q))$ is given by Rice formula,
\begin{eqnarray*}
\E[\ell({\mathcal Z}(Q))]&=&\int_Q\E\big[(\det{}^{\bot}Z(x))^{1/2}|\xi(x)=\eta(x)=0\big]p_{\xi(x),\eta(x)}(0,0)\,dx\\
&=&\vol(Q)\,\frac{1}{2\pi}\,\E\big[(\det{}^{\bot}Z(0))^{1/2}\big],
\end{eqnarray*}
where we have used stationarity and independence to get the second line as well as the fact that $\xi(0)$ and $\eta(0)$ 
are independent standard Gaussian random variables. \\
Formula \eqref{e:cross} gives
\begin{equation} \label{e:ElZQ}
\E[\ell ({\mathcal Z}(Q))]
=\frac{\vol(Q)}{2\pi}\,\E|\xi'(0)\wedge\eta'(0)|.
\end{equation}
Recall that $\Cov(\xi'_i(0),\xi'_j(0))=\Cov(\eta'_i(0),\eta'_j(0))=-r_{ij}''(0)$. 
Without loss of generality (see \cite{ad:tay}) we only study the case where
$$-r''(0)=\begin{pmatrix}\lambda_1&0&0\\0&\lambda_2&0\\0&0&\lambda_3\end{pmatrix}=D.$$ 
We write $\xi'(0)=D^{\frac12}N$ and $\eta'(0)=D^{\frac12}N'$, being $N$ and $N'$ 
two independent $N(0,I_3)$ vectors. Then, using the following algebraic property of the cross product,
\[
D^{\frac12}N\wedge D^{\frac12}N'=(\det D^{\frac12})D^{-\frac12}(N\wedge N'),
\] 
it holds
\[
\E[\ell ({\mathcal Z}(Q))]
=\vol(Q)\;\frac{\sqrt{\lambda_1\lambda_2\lambda_3}}{2\pi}\;\E|D^{-\frac12}(N\wedge N')|.
\]
This proves Theorem \ref{TheoLength}. We now move to the corollary.
\\~\\
$(i)$ If $\lambda_1=\lambda_2=\lambda_3=\lambda$, then $D=\lambda I_3$. Furthermore, from \cite{ALW11} page 34, we know that 
$\E|N\wedge N'|=2$,
thus $\E[\ell ({\mathcal Z}(Q))]=\vol(Q)\,\frac{\lambda}{\pi}.$
~\\
$(ii)$ Let $\lambda>0$ be fixed and consider $\lambda^*=(\lambda,\lambda,\lambda)$, ${\boldsymbol \lambda}=(\lambda_1,\lambda_2,\lambda_3)\in\R^3$.\\
Recall that $Z(0)=(\xi'(0),\eta'(0))\sim N(0,\diag(\lambda_1,\lambda_2,\lambda_3,\lambda_1,\lambda_2,\lambda_3))$. 
Then, from \eqref{e:ElZQ} we have
\[
\E[\ell ({\mathcal Z}(Q))]
=\frac{\vol(Q)}{2\pi}\,\int_{\R^3}\int_{\R^3}|y\wedge y'|p_{{\boldsymbol \lambda}}(y)p_{{\boldsymbol \lambda}}(y')dydy',
\]
where $p_{{\boldsymbol \lambda}}(y)=(2\pi)^{-3/2}(\lambda_1\lambda_2\lambda_3)^{-1/2}\,e^{-\frac 12\sum_{i=1}^3(\lambda_i)^{-1/2}y_i^2}$.
Hence, for $i=1,2,3$ we get
\[
\partial_{\lambda_i}\big(p_{{\boldsymbol \lambda}}(y)\big)
=\Big(-\frac{1}{2\lambda_i}+\frac{1}{4(\lambda_i)^{3/2}}y_i^2\Big)\,p_{{\boldsymbol \lambda}}(y)
\]
and so
\[
\partial_{\lambda_i}\E(|\xi'\wedge\eta'|)\big|_{{\boldsymbol \lambda}=\lambda^*}
=-\E(|N\wedge N'|)+\frac{\sqrt{\lambda}}{2}\E(|N\wedge N'|(N_i)^2)=-2+\frac{4}{3}\sqrt{\lambda},
\]
being $(N,N')$ a standard normal vector in $\R^6$.
Taylor formula allows one to terminate the proof of Corollary \ref{cor:length}.

\section{Hermite expansion and chaotic decomposition}\label{sec:H}

In this section, we introduce preliminary materials that will be useful in the sequel. 
It mainly deals with Hermite expansion which yields It\^o-Wiener's standard chaotic decomposition. 

We introduce Hermite polynomials by $H_0(x)=1$, $H_1(x)=x$ for $x\in\R$ and for $n\geq 2$ by
$$
H_n(x)=xH_{n-1}(x)-(n-1)H_{n-2}(x),\quad x\in\R.
$$
They form a complete orthogonal system in $L^2(\varphi(dx))$, 
being $\varphi$ the standard normal density function in $\R$. 
More precisely, for standard normal $X,Y$ with covariance $\rho$ 
it holds
\begin{equation}\label{e:Hiso}
\E(H_p(X)H_q(Y))=\delta_{pq}p!\rho^p,
\end{equation} 
being $\delta_{pq}$ Kronecker's delta function.
\\
The multi-dimensional Hermite polynomials are tensorial products 
of their one-dimensional versions. 
That is, for ${\boldsymbol \alpha}=(\alpha_i)_i\in\N^m$ and ${\boldsymbol y}=(y_i)_i\in\R^m$,  
$$
\tilde{H}_{\boldsymbol\alpha}({\boldsymbol y})=\prod^m_{i=1}H_{\alpha_i}(y_i).
$$
In this case, Hermite polynomials form a complete orthogonal system 
of $L^2(\varphi_m(d{\boldsymbol y}))$ 
being $\varphi_m$ the standard normal density function in $\R^m$. 
In other words, if $f\in L^2(\varphi_m(d{\boldsymbol y}))$, 
then $f$ can be written in the $L^2$-sense as
\begin{equation*}
f({\boldsymbol y})
=\sum^\infty_{q=0}\sum_{\boldsymbol \alpha \in \N^m,\,|{\boldsymbol \alpha}|=q}
f_{{\boldsymbol \alpha}}\tilde{H}_{\boldsymbol \alpha}({\boldsymbol y}),\quad 
{\boldsymbol y}\in \R^m,
\end{equation*}
with $|\balpha|=\sum^m_{i=1}\alpha_i$ and 
\begin{equation*}
f_{\boldsymbol \alpha}=\frac{1}{\balpha !}
\int_{\R^m}f({\boldsymbol y})\tilde{H}_{\boldsymbol \alpha}({\boldsymbol y})\varphi_m(d{\boldsymbol y}),
\end{equation*}
with $\balpha !=\prod^m_{i=1}\alpha_i$.\\

\medskip

We are now ready to state the Hermite expansion of the length of the zero set. From now on, we restrict our model to the isotropic case and assume that the second spectral moment $\lambda$ is positive.

Denote
\begin{equation*}
\overline{Y}(x)=\left(\xi(x),\eta(x),\frac{\xi'(x)}{\sqrt{\lambda}},\frac{\eta'(x)}{\sqrt{\lambda}}\right)\,\in\R^8. 
\end{equation*}
Let also $c_{\boldsymbol \alpha}=b_{\alpha_1}b_{\alpha_2}a_{(\alpha_3,\dots,\alpha_8)}$ being 
\begin{equation}\label{e:b-alpha}
b_\alpha=\frac{1}{\alpha!\sqrt{2\pi}}H_\alpha(0) 
\end{equation}
and $a_{(\alpha_3,\dots,\alpha_8)}$ the Hermite coefficient of ${\boldsymbol y}\in \R^6\mapsto\det{^{\bot}}({\boldsymbol y})^{1/2}$.

\begin{proposition}\label{p:exp}
With the above notations, it holds in the $L^2$-sense that
\begin{equation*}
\ell({\mathcal Z}(Q))-\E(\ell({\mathcal Z}(Q)))=\lambda\sum_{q\ge 1}I_{2q}(Q),
\end{equation*}
where 
\begin{equation*}
I_{2q}(Q)=\sum_{{\boldsymbol \alpha}\in\N^8,\\|{\boldsymbol \alpha}|=2q}
c_{\boldsymbol \alpha}\int_{Q}\tilde{H}_{\boldsymbol \alpha}(\overline{Y}(x))dx.
\end{equation*}
\end{proposition}
The proof of this proposition is based on the following standard lemma.  
\begin{lemma}
Consider a positive even kernel $h:\R\to\R$ such that $\int h=1$. 
For $\eps>0$, let $h_\eps(x)=\frac{1}{\eps}h(x/\eps)$. 
Set $\overline{h}_\eps:\R^2\to\R$ by $\overline{h}_\eps(x,y)=h_\eps(x)h_\eps(y)$. 
Define 
\begin{equation*}
\ell_\eps=\frac{1}{\eps^2}\int_{Q}\overline{h}_\eps(\xi(x),\eta(x))
(\det{}^{\bot}(Z(x)))^{1/2}dx. 
\end{equation*}
Hence, $\ell_\eps$ converge to $\ell({\mathcal Z}(Q))$ almost surely and in $L^2$. 
Besides, $\ell_\eps$ admits the Hermite ($L^2$) expansion
\begin{equation*}
\ell_\eps
=\lambda\sum^\infty_{q=1}\sum_{{\boldsymbol \alpha}\in\N^8,\\|{\boldsymbol \alpha}|=2q}
c^{\eps}_{\boldsymbol \alpha}\int_{Q}\tilde{H}_{\boldsymbol \alpha}(\overline{Y}(x))dx,
\end{equation*}
being $c^{\eps}_{\boldsymbol \alpha}=b^\eps_{\alpha_1}b^\eps_{\alpha_2}a_{(\alpha_3,\dots,\alpha_8)}$ 
with $a_{(\alpha_3,\dots,\alpha_8)}$ as above and  
$b^\eps_{\alpha}$ the Hermite coefficients of $h_\eps$.
\end{lemma}

The orthogonality of the chaotic decomposition given by Proposition \ref{p:exp} yields the following expansion for the variance of the zero set length, 
\[ 
\Var(\ell({\mathcal Z}(Q)))=\lambda^2\sum_{q\ge 1}\Var(I_{2q}(Q)).
\]
We state a lemma concerning the asymptotic behaviour of this series as $Q\uparrow \R^3$. Recall that function $R$ is defined in \eqref{defR}.

\begin{lemma} \label{L2q_0}
Let $Q_n=[-n,n]^3$. 
If the covariance function $r$ is isotropic, if $R(x)\to 0$ as $|x|\to\infty$ and 
$R$ belongs to $L^{2q_0}(\R^3)$ for some positive integer $q_0$, then there 
exists $V_{2q_0}\in [0,+\infty)$ such that
$$\lim_{n\to \infty}\,\frac{\sum_{q\ge q_0}\Var(I_{2q}(Q_n))}{\vol(Q_n)}=V_{2q_0}.$$
\end{lemma}
~\\
\begin{proof} 
For simplicity, we normalize $R$ as
\begin{equation}\label{e:Rnorm}
R(x)=\max\left\{|r(x)|,\frac{|r'_i(x)|}{\sqrt{\lambda}},\frac{|r''_{ij}(x)|}{\lambda}:1\leq i,j\leq 3\right\}.
\end{equation} 
The proof follows the same lines as that of Proposition 2.1 in \cite{EL16} with minor modifications. 

We only detail the part that needs to be adapted.  
For fixed $q\geq q_0$ we write
\begin{equation*}
\Var(I_{2q}(Q_n))=
\lambda^2\sum_{|{\boldsymbol \alpha}|=|{\boldsymbol \beta}|=2q}c_{\boldsymbol \alpha}c_{\boldsymbol \beta}
\int_{\R^3}\vol(Q_n\cap Q_n-x)\E[\tilde{H}_{\boldsymbol \alpha}(\overline{Y}(0))\tilde{H}_{\boldsymbol \beta}(\overline{Y}(x))]\,dx. 
\end{equation*}
Using Mehler's formula (see Lemma 10.7 in \cite{aw}),  we get the next upper bound for any ${\boldsymbol \alpha}$ 
and ${\boldsymbol \beta}$ in $\N^8$ such that $|{\boldsymbol \alpha}|=|{\boldsymbol 
\beta}|=2q$,
\begin{equation*}
\E[\tilde{H}_{\boldsymbol \alpha}(\overline{Y}(0))\tilde{H}_{\boldsymbol \beta}(\overline{Y}(x))]
=\sum_{\Lambda_{\balpha,\bbeta}}{\boldsymbol \alpha}!{\boldsymbol \beta}!\prod_{1\leq i,j\leq 
8}\frac{\Cov(\overline{Y}_i(0)\overline{Y}_j(x))^{d_{ij}}}{d_{ij}!}
\le K_q\,R(x)^{2q},
\end{equation*}
where 
$\Lambda_{\balpha,\bbeta}=\{d_{ij}\geq0:\sum_id_{ij}=\alpha_j,\sum_jd_{ij}=\beta_i\}$.
Here we have used that 
\[
|\Cov(\overline{Y}_i(0),\overline{Y}_j(x))|\le R(x),\quad \textrm{for any } x\in \R^3, 
\]
and that $\sum_{i,j}d_{ij}=2q$.
Thus, it follows that for any $q\ge q_0$, $\frac{\Var(I_{2q}(Q_n))}{\vol(Q_n)}$
has a finite limit as $n\to \infty$. 

The end of the proof is exactly as in \cite{EL16}. 
\end{proof}

\medskip

A key role in our asymptotic analysis of $\ell({\mathcal Z}(Q))$ will be played by the second chaotic component. 
Hence, we end this section analyzing $I_2(Q)$. 

In the next lemma, we do not assume any restrictive condition on the covariance function $r$, except it is isotropic. 

We denote by $e_j\in\N^8$ the $j$-th canonical vector, that is, the vector all of whose entries are zero but the $j$-th which is one. 
\begin{lemma}\label{l:I-2}
With previous notations and assuming $r$ is sisotropic, we have 
\[
I_2(Q)=\sum_{1\leq k\leq 8}c_{2e_k}\int_{Q}\tilde{H}_{2e_k}(\overline{Y}(x))dx,
\]
with $c_{2e_k}=-\frac{1}{2\pi}$ for $k=1,2$ and $c_{2e_k}=\frac{1}{6\pi}$ for $k=3,\dots,8$.
\\~\\
Moreover
\begin{equation}\label{e:vI2l}
 \Var(I_2(Q))=\frac{1}{\pi^2}\int_{\R^3}\vol(Q\cap Q-x){\mathcal D}r(x)dx,
\end{equation} 
where the functional ${\mathcal D}$ is defined by 
$${\mathcal D}r(x)=r(x)^2-\frac{2}{3\lambda}\sum^3_{j=1}(r'_j(x))^2+\frac{1}{9\lambda^2}\sum^3_{j,l=1}(r''_{j,l}(x))^2,\quad x\in \R^3,$$
or equivalently, writing $r(x)=\gamma(|x|)$ for some map $\gamma:\R^+\to \R$,
\begin{equation}\label{e:Dr}
{\mathcal D}r(x)=\gamma(|x|)^2+\frac{2}{3\lambda}(\frac{1}{3\lambda |x|^2}-1)\gamma'(|x|)^2+\frac{1}{9\lambda^2}\gamma''(|x|)^2,\quad x\in \R^3.
\end{equation}
\end{lemma}
~\\
\begin{proof}
From Proposition \ref{p:exp} we have
\begin{eqnarray*}
I_2(Q)&=&2\,\sum_{1\leq i< j\leq 8}c_{e_i+e_j}\int_{Q}\tilde{H}_{e_i+e_j}(\overline{Y}(x))dx
+\sum_{1\leq k\leq 8}c_{2e_k}\int_{Q}\tilde{H}_{2e_k}(\overline{Y}(x))dx\\
&:=&2\,I^{(1)}_2+I^{(2)}_2,
\end{eqnarray*}
where we recall that $c_{\boldsymbol\alpha}=b_{\alpha_1}b_{\alpha_2}a_{(\alpha_3,\dots,\alpha_8)}$. \\

Let us first show that $c_{e_i+e_j}=0$ for all $1\leq i<j\leq 8$. This will imply that $I^{(1)}_2=0$ and hence that $I_2(Q)=I^{(2)}_2$.\\
From Equation \eqref{e:b-alpha} it follows that $b_{1}=0$. Thus, $c_{e_i+e_j}=0$ for $i\leq 2$ and any $j>i$.\\
Consider $i$ and $j$ in $\{3,\dots,8\}$ with $i<j$ and $j-i\neq 3$. Then,
\begin{align*}
a_{e_i+e_j}&=\int_{\R^6}(\det{}^{\bot}{\boldsymbol y})^{1/2}H_1(y_i)H_1(y_j)\varphi_6({\boldsymbol y})d{\boldsymbol y}\\
&=\E\big(|(N_3,N_4,N_5)\wedge (N_6,N_7,N_8)|N_iN_j\big),
\end{align*}
for $N=(N_3,\dots,N_8)$ standard normal random vector in $\R^6$. 
Denote by $N'=(N'_3,\dots,N'_8)$ the vector obtained from $N$ replacing $N_i$ and $N_{i+3}$ by $-N_i$ and $-N_{i+3}$ respectively. 
It is easy to check that $|(N_3,N_4,N_5)\wedge (N_6,N_7,N_8)|=|(N'_3,N'_4,N'_5)\wedge (N'_6,N'_7,N'_8)|$. 
Since $N$ and $N'$ are equally distributed, we have
\begin{align*}
a_{e_i+e_j}&=\E\Big(|(N'_3,N'_4,N'_5)\wedge(N'_6,N'_7,N'_8)|N'_iN'_j\Big)\\
&=\E\Big(|(N_3,N_4,N_5)\wedge(N_6,N_7,N_8)|(-N_i)N_j\Big)=-a_{e_i+e_j}.
\end{align*}
Thus, $a_{e_i+e_j}=c_{e_i+e_j}=0$ if $i<j\in\{3,\dots,8\}$ with $j-i\neq 3$. 
The same argument but replacing $N$ by $N'=(-N_3,-N_4,-N_5,N_6,N_7,N_8)$ yields $c_{e_3+e_6}=c_{e_4+e_7}=c_{e_5+e_8}=0$.\\

Besides, the coefficients $c_{2e_k}$ in $I^{(2)}_2$, $k=1,\dots,8$, can be obtained by routine computations via a change to spherical coordinates.\\

Finally, we compute the variance of $I_2(Q)$. Note that in $I^{(2)}_2$, the random variables corresponding to $k\in\{1,3,4,5\}$ are independent (and equally distributed) of those corresponding to $k\in\{2,6,7,8\}$. Thus, we consider one of these two blocks.
\begin{align*}
 \Var(I_2(Q))
 &=2\sum_{j,l\in\{1,3,4,5\}}c_{2e_j}c_{2e_l}\int_{Q\times Q}\E(H_2(\overline{Y}_j(s))H_2(\overline{Y}_l(t)))dsdt \nonumber \\
 &=4\sum_{j,l\in\{1,3,4,5\}}c_{2e_j}c_{2e_l}\int_{Q\times Q}(\E\overline{Y}_j(s)\overline{Y}_l(t))^2dsdt\nonumber \\
 &=4\int_{\R^3}\vol(Q\cap Q-x)\big(\!\sum_{j,l\in\{1,3,4,5\}}c_{2e_j}c_{2e_l} (\E\overline{Y}_j(0)\overline{Y}_l(x))^2\big) dx, 
\end{align*} 
where we have used \eqref{e:Hiso} and the stationarity of $\overline{Y}(x)$. 
Since, the covariances among the coordinates of $\overline{Y}(x)$ are the corresponding derivatives of $r(x)$, 
the result follows.
\end{proof}


\section{Asymptotic variance and limit theorems} \label{sec:clt}
In this section, we estimate the asymptotic behaviour of the variance of $\ell(\mathcal Z(Q))$ as $Q\uparrow \R^3$ and derive CLT results. We consider separately the cases where $R$ (see \eqref{defR}) is square integrable on $\R^3$ or not, leading to distinct asymptotics.  

\subsection{Square integrable case} \label{sec:L2case}

We first prove Theorem \ref{TheoVar} and we postpone the exhibition of examples to the end of the section. 
\\~\\
{\em Proof of Theorem \ref{TheoVar}.}
Let us assume that the conditions of Theorem \ref{TheoVar} are satisfied. In order to simplify notations, we write $Q$ instead of $Q_n=[-n,n]^3$ and $Q\uparrow\R^3$ instead of $n\to \infty$. Note that the second spectral moment $\lambda$ does not vanish 
since it is equal to $\int_{\R^3}(k_1)^2\frac{f(\k)}{|\k|^2}d\k$, being $f$ the spectral density. \\
{\em (i)}
The upper bound for the asymptotic variance follows from Lemma \ref{L2q_0} with $q_0=1$. Thus, it remains to prove that the limit variance is strictly positive. \\
Recall that Proposition \ref{p:exp} yields
\begin{equation*}
\Var(\ell({\mathcal Z}(Q)))=\lambda^2\sum_{q\ge 1}\Var(I_{2q}(Q)) \geq \lambda^2\Var(I_2(Q)),
\end{equation*}
and that $\Var(I_{2q}(Q))$ is given by \eqref{e:vI2l} in Lemma \ref{l:I-2}.\\
Since $R\in L^2(\R^3)$ it follows that ${\mathcal D}r\in L^1(\R^3)$. Thus, by Lebesgue's dominated convergence theorem, 
\[
\lim_{Q\uparrow\R^3}\frac{ \Var(I_2(Q))}{\vol(Q)}
=\frac{1}{\pi^2}\int_{\R^3}{\mathcal D}r(x)dx.
\]
Denoting by $f$ the density of $\Pi$, Equation \eqref{e:r} now reads
\[
 r(x)=\int_{\R^3}e^{i\langle\k,x\rangle}\frac{f(\k)}{|\k|^2}d\k,\quad x\in \R^3.
\]
Taking derivatives, we get
\[
 r'_j(x)=\int_{\R^3}ik_je^{i\langle\k,x\rangle}\frac{f(\k)}{|\k|^2}d\k;\quad
 r''_{j,l}(x)=-\int_{\R^3}k_jk_le^{i\langle\k,x\rangle}\frac{f(\k)}{|\k|^2}d\k.
\]
Hence, using Plancherel identity, we get
\[ \int_{\R^3}{\mathcal D}r(x)dx=\int_{\R^3}\left(1+\frac{|\k|^2}{3\lambda}\right)^2\frac{f(\k)^2}{|\k|^4}d\k>0.
\]
Statement {\em (i)} follows. 
\\~\\
{\em (ii)} From item {\em (i)}, we know that
\[
\frac{\Var(\ell({\mathcal Z}(Q)))}{\vol(Q)}=\sum_{q\ge 0}\Var\left(\frac{I_{2q}(Q)}{\sqrt{\vol(Q)}}\right)\underset{Q\uparrow\R^3}{\longrightarrow}V_2<+\infty.
\]
Furthermore, in the same form as Proposition 2.1 in \cite{EL16}, one can prove that 
\[
\lim_{N\to\infty}\sup_{Q\subset\R^3}\sum_{q\ge N}\Var\left(\frac{I_{2q}(Q)}{\sqrt{\vol(Q)}}\right)=0.
\]
Hence, to establish the CLT for $\ell({\mathcal Z}(Q))$, 
it is sufficient to prove the asymptotic normality of each normalized component $I_{2q}(Q)/\sqrt{\vol(Q)}$ as $Q\uparrow\R^3$, 
see \cite[Th. 11.8.3]{PT11}. We do this in two steps.
\medskip

\noindent{\bf Step 1:} 
We translate the Hermite expansion obtained so far 
to the framework of isonormal processes, see \cite[Ch.8]{PT11} for the details. 

Let ${\mathcal H}={\mathcal H}_1\oplus{\mathcal H}_2$, with ${\mathcal H}_i=L^2(\R^3,\frac{\Pi(d\k)}{|\k|^2})$, $i=1,2$, endowed with the inner product
\[
 \left\langle c\oplus s,c'\oplus s'\right\rangle_{{\mathcal H}}
 =\left\langle c,c'\right\rangle_{{\mathcal H}_1}+\left\langle s,s'\right\rangle_{{\mathcal H}_2}.
\]
We also set $I^B_1:{\mathcal H}\to L^2(B)=L^2_{\R}(W)$ by 
\[
I^B_1(c\oplus s)=I^{W_1}_1(c)+I^{W_2}_1(s), 
\]
being $W_1$ and $W_2$ the real and the imaginary parts of $W$ respectively. 
It follows that
\begin{equation*}
\E\big(I^B_1(c\oplus s)\, I^B_1(c'\oplus s')\big)
=\left\langle c\oplus s,c'\oplus s'\right\rangle_{{\mathcal H}}. 
\end{equation*}
Thus $B$ is a Gaussian isonormal process. 

Now, 
let $h_{i,x}(\k)=c_{i,x}(\k)\oplus s_{i,x}(\k)\in{\mathcal H}$ be such that 
$\overline{Y}_i(x)=I^{B}_1(h_{i,x}(\k))$, $i=1,\dots,8$. 
For instance, since $\overline{Y}_1(x)=\xi(x)$, we have $h_{1,x}(\k)=c_{1,x}(\k)\oplus s_{1,x}(\k)$ with
\[
c_{1,x}(\k)=\frac{\cos(\k\cdot x)}{|\k|}~\mbox{ and }~ s_{1,x}(\k)=-\frac{\sin(\k\cdot x)}{|\k|}.
\]

Let $h_{i,x}(\k)\otimes h_{j,y}(\k')=(c_{i,x}\otimes c_{j,y}(\k,\k'))\oplus(s_{i,x}\otimes s_{j,y}(\k,\k'))$. 
By definition of the $2q$-folded multiple Wiener integral with respect to $B$, we get
\begin{equation*}
\tilde{H}_{{\boldsymbol\alpha}}(\overline{Y}(x))
=\prod^{8}_{i=1}H_{\alpha_i}(\overline{Y}_i(x))
=I^B_{2q}\left(\otimes^{8}_{i=1}h_{i,x}^{\alpha_i}\right),
\end{equation*}
where $|\balpha|=2q$ and 
$\otimes^{8}_{i=1}h_{i,x}^{\alpha_i}=\otimes^{8}_{i=1}h_{i,x}^{\alpha_i}(\K)$ 
stands for the tensorial products of the kernels $h_{i,x}$ 
for $\K=(\k_1,\dots,\k_{2q})\in(\R^3)^{2q}$. \\
Therefore,
\begin{equation*}
I_{2q}(Q)=I^{B}_{2q}\left(g_{2q}\right), 
\end{equation*}
with
\[
g_{2q}(\K)=\sum_{{\boldsymbol \alpha}\in\N^8,|{\boldsymbol \alpha}|=2q}c_{\boldsymbol \alpha}\int_{Q}
\otimes^{8}_{i=1}h_{i,x}^{\alpha_i}(\K)dx.
\]
\medskip

\noindent{\bf Step 2:} 
Once that $I_{2q}(Q)$ has been written as a multiple integral, 
thanks to the fourth moment Theorem \cite[Th. 6.3.1]{NP12}, 
to establish its asymptotic normality, it suffices to prove that the $2$-norms 
of the so-called contractions of the normalized kernels 
tend to $0$. 

\begin{remark}
 There are other ways of proving the asymptotic normality 
 of a sequence of random variables living in a fixed chaos,  
 see \cite{NP12,PT11} for details. 
 We choose contractions since the computations are straightforward 
 in the present case.
\end{remark}

Let us recall that for $p\in\N$, symmetric $f,g\in {\mathcal H}^{\otimes p}$ and  $1\leq n\leq p$, 
the $n$-th contraction is defined as
\begin{equation*}
f\otimes_n g
=\sum_{i_1,\dots,i_n=1}^{\infty}\<f,e_{i_1}\otimes\dots\otimes e_{i_n}\>_{{\mathcal H}^{\otimes n}}\otimes
\<g,e_{i_1}\otimes\dots\otimes e_{i_n}\>_{{\mathcal H}^{\otimes n}},
\end{equation*}
being $\{e_i\}_i$ a complete orthogonal system in ${\mathcal H}$. 
The definition does not depend on the choice of the basis $\{e_i\}_i$.
Since each covariance is bounded by $R$, in order to avoid messy notations 
we do not symmetrize the kernels in the next lines.

Note that in the case that $f=\otimes^{p}_{i=1}f_i$ and $g=\otimes^{p}_{i=1}g_i$, 
then
\begin{equation}\label{e:contr}
f\otimes_n g=\prod^{n}_{i=1}\left\langle f_i,g_i\right\rangle_{\mathcal H}\left(
\otimes^{p-n}_{i=1}f_i\otimes\otimes^{2p-2n}_{i=p-n+1}g_i\right).
\end{equation}
\medskip

In our case, $p=2q$ and
\[
g_{2q}\otimes_{n}g_{2q}
=\sum_{|\balpha|=|\balpha'|=2q}c_{\balpha}c_{\balpha'}\int_{Q\times Q}
\big(\otimes^{8}_{i=1}h^{\alpha_i}_{i,x}\otimes_n\otimes^{8}_{i=1}h^{\alpha'_i}_{i,x'}\big)dxdx'.
\]
Besides, from \eqref{e:contr} we see that 
the contraction in the last integral 
yields $n$ inner products (using $n$ kernels with $x$ and $n$ kernels with $x'$)
that, 
since $I^B_1$ is an isonormal process, 
equal the covariances of the corresponding 
elements of $\overline{Y}(x)$ and $\overline{Y}(x')$. 
For instance, 
$\left\langle h_{1,x},h_{1,x'}\right\rangle_{{\mathcal H}}=\E(\xi(x)\xi(x'))=r(x-x')$. 
and $\left\langle h_{1,x},h_{3,x'}\right\rangle_{{\mathcal H}}=\E(\xi(x)\xi'_1(x'))=r'_1(x-x')$, etc. 
Furthermore, it remains 'un-used' $2q-n$ kernels of $x$ and $2q-n$ of $x'$.\\

Recall that $R(x)\to0$ as $|x|\to\infty$ and that $R\in L^{2}(\R^3)$. 

Taking ${\mathcal H}^{4q-2n}$ norms and using the fact that 
all the covariances of $\overline{Y}$ are bounded by $R$, we get
\begin{multline*}
\left\|\frac{g_{2q}}{\sqrt{\vol(Q)}}\otimes_n \frac{g_{2q}}{\sqrt{\vol(Q)}}\right\|^2\\
\leq \frac{C_q}{\vol(Q)^2}\int_{Q^4}
R^n(x-x')R^n(y-y')R^{2q-n}(x-y)R^{2q-n}(x'-y')dxdx'dydy',
\end{multline*}
where $C_q$ is some constant which takes into account the coefficients $c_{\balpha}$ and the number of terms in the sums.

Now, we make the isometric change of variables 
$(x,x',y,y')\mapsto (x-x',y-y',x-y,x')$. 
Next, we enlarge the domain of integration to $\widetilde{Q}^4$ 
so that it includes the image of $Q^4$ under the change of variables
and $\vol(\widetilde{Q})=c\vol(Q)$ for some constant $c$. 
Hence, we get
\begin{multline*}
\left\|\frac{g_{2q}}{\sqrt{\vol(Q)}}\otimes_n \frac{g_{2q}}{\sqrt{\vol(Q)}}\right\|^2\\
\leq \frac{C_q}{\vol(Q)^2}\int_{\widetilde{Q}^4}
R^{n}(u_1)R^{n}(u_2)R^{2q-n}(u_3)R^{2q-n}(u_2+u_3-u_1)du_1du_2du_3du_4\\
\leq \frac{C_q}{\vol(Q)}\int_{\widetilde{Q}^3}
R^n(u_1)R^n(u_2)R^{2q-n}(u_3)R^{2q-n}(u_2+u_3-u_1)du_1du_2du_3.
\end{multline*}
If $1<n<2q-1$ (thus $q>1$), 
since $R\in L^2$ it follows that  the contractions tend to $0$.\\ 
Now, assume that $n=1$ and $q=1$ which is the most difficult case. 
By Cauchy-Schwarz, for fixed $u_3$ and $u_1$, $\int_{\widetilde{Q}}R(u_2)R(u_2+u_3-u_1)du_2$ is bounded. 
Hence, it suffices to prove that as $Q\uparrow\R^3$
\begin{equation}\label{e:L2c}
\frac1{\sqrt{\vol(Q)}}\int_{Q}R(u)du\to0. 
\end{equation}
To see this, take $Q'_n\subset Q_n$ such that $Q'_n\uparrow\R^3$ with $\vol(Q'_n)=o(\sqrt{\vol(Q_n)})$. 
Thus
\begin{align*}
\frac1{\sqrt{\vol(Q_n)}}\,\int_{Q_n}R(u)du
&=\frac1{\sqrt{\vol(Q_n)}}\,\int_{Q_n\setminus Q'_n}R(u)du+o(1)\\
&=\frac{\vol(Q_n\setminus Q'_n)}{\sqrt{\vol(Q_n)}}\,\int_{Q_n\setminus Q'_n}R(u)\frac{du}{\vol(Q_n\setminus Q'_n)}+o(1)\\
&\leq \sqrt{\frac{\vol(Q_n\setminus Q'_n)}{\vol(Q_n)}}\,\left[\int_{Q_n\setminus Q'_n}R^2(u)du\right]^{1/2}+o(1)\\  
&\leq  \left[\int_{(Q'_n)^c}R^2(u)du\right]^{1/2}+o(1)
\to_n0,
\end{align*}
where the first inequality is due to Jensen's inequality. Hence, \eqref{e:L2c} follows. 
The remaining cases are similar and easier. 

Hence,
\[
\left\|\frac{g_{2q}}{\sqrt{\vol(Q)}}\otimes_n \frac{g_{2q}}{\sqrt{\vol(Q)}}\right\|^2
 \to_{Q\uparrow\R^3}0.
\]
This completes the proof of the CLT assertion in Theorem \ref{TheoVar}. 
\hfill $\Box$
\\~\\
In order to illustrate the results, we end this section by giving three examples of random waves models that enter in the square integrable case.
We use formulas \eqref{e:r} and \eqref{e:riso} to express the covariance function $r$ in terms of the power spectrum $\Pi$. 
In all our examples we assume that the power spectrum admits a density with respect to Lebesgue measure 
and since we focus on isotropic examples we write it as $f(|.|)$. Hence the two next identities will be in force
\[r(x)=\int_{\R^3}\exp(i\langle \mathbf{k},x\rangle)\,\frac{f(|\k|)}{|\k|^2}\,d\k
=4\pi\,\int_{\R^+}\frac{\sin(\rho|x|)}{\rho|x|}\,f(\rho)\,d\rho,\]
with normalization $\int_{\R^3}\frac{1}{|\k|^2}f(|\k|)\,d\k=4\pi\,\int_{\R^+}f(\rho)\,d\rho=1$. 
\\~\\
{\bf Bargmann-Fock model.} Let us take $f(\rho)=(2\pi)^{-3/2}\,\rho^2\,e^{-\rho^2/2},~\rho\in \R^+$ as spectral density. In this case,
\[r(x)=(2\pi)^{-3/2}\,\int_{\R^3}\exp(i\k.x)\,e^{-|\k|^2/2}\,d\k=e^{-|x|^2/2},~x\in \R^3.\]
Since the covariance function as well as all its derivatives belong to all $L^p(\R^3)$, Theorem \ref{TheoVar} applies. 
\\~\\
{\bf Gamma type.} Let us take $f(\rho)=\frac{\beta^{p+1}}{4\pi p!}\,\rho^p\,e^{-\beta\rho},~\rho\in \R^+$ with $p$ a positive integer and 
$\beta$ some positive real constant. We remark that $f(\rho)=\frac{\beta}{4\pi p}\rho\gamma(\rho),$ where $\gamma$ is the probability density 
function 
of a $\Gamma(p,\beta)$-distribution. We then write the covariance function as  
\[
r(x)=\frac{\beta}{p|x|}\,\int_{\R^+}\sin(\rho|x|)\,\gamma(\rho)\,d\rho=\frac{\beta}{p|x|}\,Im\big(\widehat{\gamma}(|x|) \big),\]
where $Im$ stands for the imaginary part of any complex number and $\widehat{\gamma}$ stands for the characteristic function of distribution 
$\gamma$. 
Since $\widehat{\gamma}(t)=(1-i\frac{|t|}{\beta})^{-p}$, we get
\begin{align*}
r(x)&= 
\frac 1p(1+\frac{|x|^2}{\beta^2})^{-p}\,\sum_{1\le j\le p;\,j\,odd}(-1)^{(j-1)/2}\,\binom{p}{j}\,\beta^{-(j-1)}\,|x|^{j-1}.  
\end{align*}
Concerning integrability properties of $r$, we note that as $|x|\to \infty$,
$$|r(x)|\approx |x|^{-(p+1)}~\mbox{ if $p$ is odd }~;~\approx|x|^{-(p+2)}~\mbox{ if $p$ is even},$$
where we denote $f(x)\approx g(x)$ for the existence of a positive constant $c$ such that $\lim_{|x|\to \infty} \frac{f(x)}{g(x)}=c$. 
In the same vein,  for odd $p$, $|r'(x)|\approx|x|^{-(p+2)}$ and $|r''(x)|\approx |x|^{-(p+3)}$ whereas for even $p$,
$|r'(x)|\approx |x|^{-(p+3)}$ and $|r''(x)|\approx |x|^{-(p+4)}.$ Hence, for $p\ge 1$, it is clear that $r$ and its derivatives belong to $L^2(\R^3)$ 
and Theorem \ref{TheoVar} again applies.
\\~\\
{\bf Black-Body radiation.}
The {\em Black-Body model} is prescribed by $f(\rho)=\tfrac{c\rho^3}{e^{\rho}-1}$, being $c$ a convenient constant. 
According to Equation (6.8) in \cite{BD00}, see also formula 2 in section 3.911 \cite{GR},
\[
r(x)=\frac{c_1}{|x|^2}-\frac{c_2 |x|\cosh( |x|)}{\sinh( |x|)^2}.
\]
This implies that $r$ and its derivatives are in $L^2(\R^3)$ and Theorem \ref{TheoVar} once more applies.\\

\medskip

In the next two subsections, we focus on two examples of random waves that behave in very different ways than the previous examples.

\subsection{Berry's monochromatic random waves model.} \label{sec:Berry}

Berry's monochromatic random waves model is defined as in \eqref{e:psi} with the power spectrum $\Pi$ 
that is uniformely distributed on the two-dimensional sphere $\S^2$.
For this isotropic model, relation \eqref{eq:Pirad} holds with $\Pi^{rad}$ proportional to the Dirac mass at 1, 
{\em i.e.} $\Pi^{rad}=\frac{1}{4\pi}\delta_1$. Thus, the covariance function is given by 
$$r(x)=\sinc(|x|),\quad x\in \R^3.$$
In particular, we stress that $r$, and hence $R$, is not square integrable on $\R^3$ and Theorem \ref{TheoVar} does not apply. Nevertheless a similar CLT holds, as stated in the next proposition.

\begin{proposition} \label{TheoMono}
Let $\psi$ be the isotropic Berry's monochromatic random wave.
Assume also that $Q_n=[-n,n]^3$. Then,
\begin{enumerate}
\item[(i)] $\displaystyle \lim_{n\to \infty}\,\frac{\Var(I_{2}(Q_n))}{\vol(Q_n)}=0$
\item[(ii)] there exists $V\in [0,+\infty)$ such that as $n\to\infty$, the distribution of  
\[ \frac{\ell({\mathcal{Z}(Q_n)})-\E(\ell({\mathcal{Z}(Q_n)}))}{\vol(Q_n)^{1/2}} \]
converges towards the centered normal distribution with variance $V$.
\end{enumerate}
\end{proposition}

Item $(ii)$ in Proposition \ref{TheoMono} states that the variance of the nodal length on a domain $Q\subset \R^3$ grows up to infinity with the same order of magnitude as the volume of $Q$. Let us mention that this behaviour strongly differs from the two-dimensional case where the variance of the nodal length on a domain $Q\subset \R^2$ is asymptotically proportional to $area(Q)\,\log(area(Q))$ as $Q$ grows up to $\R^2$ (see \cite{Be02,NPR17}). However, the vanishing second chaotic component that is observed in 2D still holds in 3D as showed by item $(i)$.

We also note that Proposition \ref{TheoMono} can be translated in terms of {\em high energy} asymptotics as in Remark \ref{rem:highNRJ} by considering Berry's monochromatic random waves with covariance $\sinc(\kappa |\cdot|)$ on the fixed domain $[-1,1]^3$ and letting $\kappa$ go to infinity. 
\\~\\
\begin{proof} 
$(i)$ We write  $Q_n=\{nx\,:\,x\in Q_1\}$ and
\[\frac{\vol(Q_n\cap (Q_n-x))}{\vol(Q_n)}=\frac{\vol(Q_1\cap (Q_1-n^{-1}x))}{\vol(Q_1)}=c(n^{-1}x),\]
where $c:y\in \R^3\mapsto c(y):=\frac{\vol(Q_1\cap Q_1-y)}{\vol(Q_1)}$ is continuous and compactly supported. 
Then, on the one hand, by Lemma \ref{l:I-2}
\begin{align}\label{e:VI2Qn}
\frac{\Var(I_2(Q_n))}{\vol(Q_n)}
&=\frac{1}{\pi^2}\int_{\R^3}c(n^{-1}x){\mathcal D}r(x)dx=\frac{4}{\pi}\int_{\R^+}C(n^{-1}\rho)D(\rho)\rho^2d\rho, 
\end{align}
where we have changed to polar coordinates and have set ${\mathcal D}r(x)=D(|x|)$ and $C(\rho)=\frac1{4\pi}\int_{\S^2}c(\rho u)d\sigma(u)$. Let us remark that $C$ is compactly supported and that $C(0)=1$.\\
On the other hand, since $\lambda=1/3$ in that case, 
one can write from \eqref{e:Dr} 
\begin{equation} \label{e:D}
D(y)\,y^2=-2\cos(2y)+4\frac{\sin(2y)}{y}+6F(y),
\end{equation}
where 
$$F(y)=\frac{1}{y^2}(\cos(2y)-\frac{\sin(2y)}{y}+\frac{\sin^2(y)}{y^2})$$
is an integrable function on $\R^+$. We now use \eqref{e:D} to split the integral in r.h.s. of \eqref{e:VI2Qn} into three terms: 
\begin{itemize}
\item Integrating twice by parts the first term yields 
$$-2\int_{\R^+}C(n^{-1}y)\,\cos(2y)\,dy=\frac{1}{n}\big(C'(0)+\int_{\R^+}\cos(2ny)\,C''(y)\,dy\big)\underset{n\to \infty}{\longrightarrow}0,$$
where we have used that $C'$ and $C''$ are compactly supported. 
\item For the second term, writing $\frac{2\sin(y)}{y}$ as the Fourier transform of the indicator function of $[-1,1]$ and using Parseval identity, one can prove that 
$$4\int_{\R^+}C(n^{-1}y)\frac{\sin(2y)}{y}dy \underset{n\to \infty}{\longrightarrow}4C(0)\frac{\pi}{2}=2\pi.$$
\item We use Lebesgue dominated convergence theorem to get the limit of the last term as $n$ goes to $\infty$: 
$$\int_{\R^+}C(n^{-1}y)\,F(y)\,dy \to \int_{\R^+}\,F(y)\,dy:=J,$$
where a tricky integration by part allows one to get that $J=-\frac{\pi}{3}$.
\end{itemize}
Finally, we conclude that 
$$\int_{\R^{+}}C(n^{-1}\rho) D(\rho)\rho^2d\rho \underset{n\to \infty}{\longrightarrow} 0+2\pi-6\frac{\pi}{3}=0.$$
and hence Part $(i)$ of Proposition \ref{TheoMono}, is now established. 
\\~\\
$(ii)$ Let us remark that $R(x)$ behaves like $1/|x|$ as $|x|\to \infty$, so that $R(x)\to 0$ and $R$ belongs to $L^4(\R^3)$. Hence, thanks to Lemma \ref{L2q_0}, we get 
\[
\lim_{n\to \infty}\,\frac{\sum_{q\ge 2}\Var(I_{2q}(Q_n))}{\vol(Q_n)}=V_4 \in [0,+\infty). 
\]
Since $\Var(\ell(\mathcal Z(Q))= \sum_{q\ge 1}\Var(I_{2q}(Q))$, applying $(i)$, we get that $\frac{\Var\ell(\mathcal Z(Q_n))}{\vol(Q_n)}\to V_4<+\infty$. 
In order to prove the CLT result, we use a similar procedure as for the proof of item $(ii)$ of Theorem \ref{TheoVar}. The difference relies on the fact that the second component 
$I_2(Q)$ in the chaotic expansion of $\ell({\mathcal Z}(Q))$ is now negligible with respect to $\sqrt{\vol(Q)}$, so we 
must only consider the contractions $\frac{g_{2q}}{\sqrt{\vol(Q)}}\otimes_n \frac{g_{2q}}{\sqrt{\vol(Q)}}$ as above for $q>1$. Since $R$ belongs to $L^4(\R^3)$, the same arguments allow us to conclude.
\end{proof}

\medskip

\subsection{Power law model.} \label{sec:power}
Our last example is a {\em power law} model named after the spectrum density given by $f(\rho)=\frac{1-\beta}{4\pi}\rho^{-\beta}\,\ind_{(0,1)}(\rho)$ with $0<\beta<1$. Using a change of variable provides the covariance function of this model as follows, 
$$r(x)=(1-\beta)\,|x|^{\beta-1}\,\int_0^{|x|}\rho^{-\beta-1}\,\sin \rho\,d\rho,\quad x\in \R^3.$$
Since the integral has a finite limit as $|x|$ tends to infinity, we get that $r(x)\approx |x|^{\beta-1}$. Hence, $r\notin L^2(\R^3)$ and one cannot apply Theorem \ref{TheoVar}. Nevertheless, for $0<\beta<1/4$ an asymptotic behaviour can be established as stated in the next proposition. 

\begin{proposition} \label{prop:power}
Let $\psi$ be a power law random waves model with parameter $\beta \in (0,1/4)$.
Assume also that $Q_n=[-n,n]^3$. Then,
\begin{enumerate}
\item[(i)] $\displaystyle \lim_{n\to \infty}\,\frac{\Var(I_{2}(Q_n))}{\vol(Q_n)^{(2\beta+4)/3}}=V \in (0,+\infty)$
\item[(ii)] as $n\to\infty$, $\displaystyle  
\frac{\ell({\mathcal{Z}(Q_n)})-\E(\ell({\mathcal{Z}(Q_n)}))}{\vol(Q_n)^{(\beta+2)/3}}$
converges in distribution towards a Rosenblatt process given by the dobble Wiener integral \eqref{eq:Ros} below.
\end{enumerate}
\end{proposition}
~\\

Note the unusual normalizing power of $\vol(Q_n)$ in the first item of Proposition \ref{prop:power}. Note also that a non-Gaussian limit is appearing in the second item, which is in hard constrast with the preceeding examples.

While proving Proposition \ref{prop:power}, we will show that  $\ell({\mathcal{Z}(Q_n)})-\E(\ell({\mathcal{Z}(Q_n)}))$ behaves as $n^{2+\beta}I_2(Q_n)$ as $n\to \infty$. This asymptotics yields the predominance of the second chaos and explains the non-Gaussian distribution limit of the normalized length. More precisely, the limit distribution belongs to the second Wiener chaos and can be written as a Rosenblatt process as introduced by Taqqu in \cite{Taq} for Hermite processes of rank two.\\

\begin{proof}
$(i)$ Since $r(x)\approx |x|^{\beta-1}$, we have that $\mathcal 
Dr(x)\approx |x|^{2\beta-2}$. So, we get for $B(0,n)$ the Euclidean ball in $\R^3$,
$\int_{B(0,n)}\mathcal Dr(x)dx\approx n^{2\beta+1}$ and hence Lemma \ref{l:I-2} yields 
\[ \Var(I_2(B(0,n) \approx (vol(B(0,n))^{(2\beta+4)/3},~n\to +\infty.\]
Replacing the ball $B(0,n)$ by the rectangle $[-n,n]^3$ does not change the order of magnitude. 
~\\
$(ii)$ 
We now deal with the asymptotic distribution. 

On the one hand, since $0<\beta<1/4$, one has $R\in L^4(\R^3)$ and Lemma \ref{L2q_0} does apply with $q_0=2$. 
Then, $\frac{\Var(\sum_{q\ge 2}I_{2q}(Q_n))}{\vol(Q_n)^{(2\beta+4)/3}}$ tends to 0 and hence, in view of the distribution limit of the normalized length, only the second chaotic component is relevant. 

On the other hand, 
by Lemma \ref{l:I-2}, $I_2(Q_n)$ is equal to the sum of two independent random variables with the same distribution.  
So we only consider one of these terms, namely $\sum_{k=1,3,4,5}c_{2e_k}\int_{Q_n}\tilde{H}_{2e_k}(\overline Y(x))dx$. \\
Thus, the first addend is constructed by using
$$\xi(x)=\int_{\R^3}e^{i<x,\mathbf \k>}\sqrt{f(|\k|)}\frac1{|\k|}dW(\mathbf \k),$$ 
being $W$ a standard complex Brownian noise. 
In particular
$$H_2(\xi(x))=\int_{\R^3\times\R^3}e^{i<x,\k+\k'>} 
\sqrt{f(|\k|)f(|\k'|)}\frac1{|\k|}\frac1{|\k'|}dW(\k)dW(\k').$$
Considering the derivatives of $\xi$, which can be written as 
$$\xi'_j(x)=i\int_{\R^3}e^{i<x,\mathbf \k>}\k_{j}\sqrt{f(|\k|)}\frac1{|\k|}dW(\mathbf \k),~j=1,2,3,$$ 
we get
$$H_2(\frac{\xi'_j(x)}{\sqrt{\lambda}})=-\frac{1}{\lambda}
\int_{\R^3\times\R^3}e^{i<x,\k+\k'>}\k_{j}\k'_{j}\sqrt{f(|\k|)f(|\k'|)}\frac1{|\k|}\frac1{|\k'|}dW(\k)dW(\k').$$
Introducing the notation 
$$g(\k,\k')= -\frac{1}{2\pi}(1+\frac{1}{3\lambda}\sum_{j=1}^3\k_{j}\k'_{j})\sqrt{f(|\k|)f(|\k'|)}\frac1{|\k||\k'|},$$
the term of our interest is
\begin{align*}
I_2(Q_n)&=\int_{Q_n}\int_{\R^3\times\R^3}e^{i<x,\k+\k'>}\,g(\k,\k')\,dW(\k)dW(\k')dx\\
&=\int_{\R^3\times\R^3}\big(\int_{Q_n}e^{i<x,\k+\k'>}dx\big)\,g(\k,\k')\,dW(\k)dW(\k')\\
&=\int_{\R^3\times\R^3}8n^3\prod_{j=1}^3\sinc(n(\k_j+\k'_j))\,g(\k,\k')\,dW(\k)dW(\k')\\
&\stackrel{d}=\int_{\R^3\times\R^3}8\prod_{j=1}^3\sinc(\k_j+\k'_j)\,g(\frac{\k}{n},\frac{\k'}{n})\,dW(\k)dW(\k'),
\end{align*}
where the change of variable $(\k,\k')\to (n\k,n\k')$ as well as the usual scaling property for Brownian measure allowed us to obtain the last identity. \\
Then, keeping in mind that $f(\rho)=\frac{1-\beta}{4\pi}\rho^{-\beta}\,\ind_{(0,1)}(\rho)$, we have
$$n^{-(2+\beta)}\,g(\frac{\k}{n},\frac{\k'}{n})\underset{n\to\infty}{\to}-\frac{1-\beta}{8\pi^2}(|\k||\k'|)^{-1-\beta/2}.$$
Hence, Theorem 1' of Dobrushin \& Major \cite{DM} yields the convergence in distribution of $n^{-(2+\beta)}I_2(Q_n)$ towards 
\begin{equation} \label{eq:Ros}
-\frac{1-\beta}{\pi^2}\int_{\R^3\times\R^3}\prod_{j=1}^3\sinc(\k_j+\k'_j)(|\k||\k'|)^{-1-\beta/2}dW(\k)dW(\k').
\end{equation}
\end{proof}


\end{document}